\documentclass[letterpaper, 10 pt, conference]{ieeeconf}  

\IEEEoverridecommandlockouts                              
\newtheorem{innercustomass}{Assumption}
\newenvironment{customass}[1]
  {\renewcommand\theinnercustomass{#1}\innercustomass}
  {\endinnercustomass}

\usepackage{graphics} 
\usepackage{epsfig} 
\usepackage{cite}
\usepackage{comment}
\graphicspath{{figures/}}

\usepackage{amsmath} 
\usepackage{amssymb}  

\usepackage{mathrsfs}

\title{\LARGE \bf
Energy-Optimal Coordination of Connected and Automated Vehicles at Multiple Intersections}

\author{A M Ishtiaque Mahbub, {\itshape{Student Member, IEEE}}, Liuhui Zhao, {\itshape{Member, IEEE}}, Dimitris Assanis, {\itshape{Member, IEEE}}, \\Andreas A. Malikopoulos, {\itshape{Senior Member, IEEE}} 
\thanks{This research was supported by ARPAE's NEXTCAR program under the award number DE-AR0000796.}
\thanks{The authors are with the Department of Mechanical Engineering, University of Delaware, Newark, DE 19716 USA (emails: \texttt{mahbub@udel.edu}; \texttt{lhzhao@udel.edu}; \texttt{dimitris@udel.edu};
	\texttt{andreas@udel.edu}.)}}

\begin{document}

\maketitle
\thispagestyle{empty}

\begin{abstract}
Urban intersections, merging roadways, roundabouts, and speed reduction zones along with the driver responses to various disturbances are the primary sources of bottlenecks in corridors that contribute to traffic congestion. The implementation of connected and automated technologies can enable a novel computational framework for real-time control aimed at optimizing energy consumption and travel time. In this paper, we propose a decentralized energy-efficient optimal control framework for two adjacent intersections. We derive a closed-form analytical solution that includes interior boundary conditions and evaluate the effectiveness of the solution through simulation. Fuel consumption and travel time are significantly reduced compared to the baseline scenario designed with conventional fixed time signalized intersections.
\end{abstract}

\indent


\section{Introduction} \label{sec:1}
The increasing traffic volume in urban areas has reached the capacity of current infrastructure resulting in congestion. Fuel efficiency and travel time can also be seriously affected in daily commute \cite{Schrank2015}. The interconnectivity of mobility systems with connected automated vehicles (CAVs) enables a novel computational framework to process massive amount of data and deliver real-time control actions that optimize energy consumption and associated benefits. CAVs can alleviate congestion at the major transportation segments such as urban intersections, merging roadways, roundabouts, and speed reduction zones, which are the primary sources of bottlenecks that contribute to traffic congestion \cite{Malikopoulos2013}.
 
Use of traffic lights is the most conventional and prevalent method for controlling traffic flow through intersections. Several research efforts have been reported in the literature proposing different approaches on coordinating CAVs through urban intersection including fuzzy logic \cite{Choi2002}, genetic algorithms \cite{Liu2007}, reservation scheme \cite{Dresner2008,DeLaFortelle2010, Huang2012}, vehicle coordination scheme \cite{Zohdy2012,Yan2009,kim2014} and swarm optimization algorithms \cite{Dong2006}. A detailed discussion of these research efforts reported in the literature to date can be found in \cite{Malikopoulos2016a}.

In earlier work, a decentralized optimal control framework was established for coordinating online CAVs in different transportation segments. A closed-form, analytical solution without considering state and control constraints was presented in \cite{Rios-Torres2015}, \cite{Rios-Torres2}, and \cite{Ntousakis:2016aa} for coordinating CAVs at highway on-ramps, in \cite{Malikopoulos2018c} at speed reduction zone, in \cite{Zhang2016a,Malikopoulos2017} at intersections, and in \cite{Zhao2018} at roundabouts. These approaches have, however, focused only on a single conflict zone, without addressing the problem of optimizing a corridor including more than one such scenarios. Thus, the development of a closed form analytical solution for a corridor consisting of multiple intersections with a single global coordinator has remained an open problem.

In this paper, we address the problem of coordinating CAVs at two adjacent intersections by 1) developing a vehicle coordination policy for throughput maximization in the corridor without any traffic lights, and 2) deriving a closed form analytical solution to the energy-optimization problem considering the interior boundary constraints that yields the optimal control input for each CAV.

The remainder of the paper is organized as follows. In Section II, we formulate the problem and provide the modeling framework. In Section III, we derive the analytical, closed form solution for the corridor control with interior constraints. In Section IV, we validate the effectiveness of the efficiency of the analytical solution in a simulation environment. Finally, the concluding remarks and discussion are provided in Section V.

\section{Problem Formulation} \label{sec:2}

\begin{figure}[ht]
\centering
\includegraphics[scale=0.5]{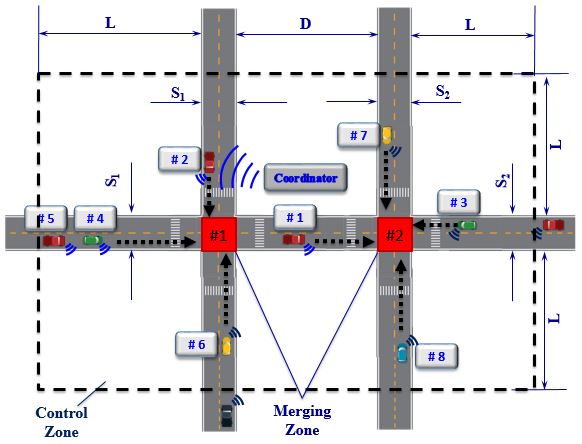} \caption{Corridor with
connected and automated vehicles.}%
\label{fig:1}%
\end{figure}

We consider a corridor (Fig. \ref{fig:1}) consisting of two adjacent urban intersections separated by a length $D$. Each intersection includes an area of potential lateral collision defined as the \textit{merging zone}, shown by the red squares of length $S_z$ for merging zone $z$, $z=1,2$, in Fig. \ref{fig:1}. The length and geometry of the merging zones are not restrictive. Both intersections are located within a \textit{control zone} illustrated in Fig. \ref{fig:1}, inside of which the CAVs can communicate with each other and with a coordinator. The distance between the entry of the control zone and the entry of the merging zone $z$ is denoted by $L_z$. Thus, the distance from the entry of the control zone to the nearest and farthest entry of the  merging zone is $L_z = L$ and $L_z=L+S_{z}+D$ respectively (Fig. \ref{fig:1}). When a CAV enters the control zone, it  exchanges information with other CAVs as well as the coordinator to derive its optimal control input (acceleration/deceleration) to cross the intersections without any rear-end or lateral collision. Note that, the coordinator only facilitates the communication among the CAVs and is not involved in any decision making process. 

In our framework, we impose the following assumptions:


\begin{customass}{1}	\label{ass:no-turn}
	We do not consider left/right turns or lane changes inside the control zone.
\end{customass}

\begin{customass}{2}	\label{ass:communication}
	Communication among CAVs occurs without any delays and errors. Each CAV $i$ is equipped with sensors to measure and share their local information.
\end{customass}

The first assumption is imposed to simplify the problem and focus on the implications of the analytical solution without adding more degrees of complexity. The second assumption may be strong, but it is relatively straightforward to relax it as long as the noise in the measurements and/or delays is bounded.

Let $N(t)\in\mathbb{N}$ be the number of CAVs inside the control zone of the corridor at time $t\in\mathbb{R}^{+}$ and $z=1,2$ be the number of merging zones of the intersections. When a CAV enters the control zone of the corridor, the coordinator receives its route information and assigns a unique identification number $i\in\mathbb{N}$. Let $\mathcal{N}(t)=\{1,..., N(t)\}$ be the merging sequence that the vehicles enter the control zone. 

\textit{Definition 1:}
\label{def:Q}
$\mathcal{Q}^z(t)$ is a local merging sequence such that $\mathcal{Q}^z(t)\subset \mathcal{N}(t)$ and includes all CAVs $M^z(t)\in \mathbb{N}$ at time $t\in \mathbb{R}^+$ that will be entering the merging zone $z$, $z={1,2}$. 

For example, CAV $\# 2,3,4,5,6\in \mathcal{Q}^1(t)$ but CAV $\# 1,7,8 \notin \mathcal{Q}^1(t)$ (Fig. \ref{fig:1}). Similarly, CAV $\# 1, 3, 4, 5, 7, 8\in \mathcal{Q}^2(t)$ but CAV $\# 2, 6 \notin \mathcal{Q}^2(t)$ (Fig. \ref{fig:1}).
If a CAV $i$ enters the control zone at time $t_i^0$ with a route designating to cross the merging zone $z$, i.e., $i\in \mathcal{Q}^z(t)$, it needs to compute the optimal time $t_i^{m_z^*}$ to enter the merging zone $z,z=1,2,$ to avoid collision.

The order of CAV $i \in \mathcal{Q}^z(t)$ satisfies the following condition,
\begin{equation}\label{eq:MS}
t_{i}^{m_z^*} \ge t_{i-1}^{m_z^*}, ~ \forall i\in \mathcal{Q}^z(t), \, i>1.
\end{equation}

When a CAV $ i\in \mathcal{Q}^z(t)$ enters the control zone at time $t_i^0$, the merging time $t_i^{m_z} = t_i^{0}+\frac{L_z}{v_i(t_i^0)} $ corresponding to its initial constant speed is compared to the optimal merging time $t_{i-1}^{m_z^*}$ of the previous CAV in the queue. If the following condition holds

\begin{equation} \label{eq:MS_check}
	t_{i}^{m_z} \ge t_{i-1}^{m_z^*}, ~ i\in \mathcal{Q}^z(t), \, i>1,
\end{equation}
then the merging sequence $\mathcal{Q}^z(t)$ is remained unchanged. However, if \eqref{eq:MS_check} does not hold, $\mathcal{Q}^z(t)$ is updated to change the order of the CAVs entering the merging zone so that \eqref{eq:MS} is not violated. The policy through which the merging sequence $\mathcal{Q}^z(t)$ is specified can be obtained as the result of an upper level vehicle coordination problem, as described in the following section. When the merging time $t_{i}^{m_z}$ is fixed, each vehicle solves a lower-level energy minimization problem that yields an analytical, closed-form optimal solution.

\subsection{Modeling Framework and Constraints}
Let $t_i^0$ be the initial time that vehicle $i$ enters the control zone of the corridor, $t_{i}^{m_z}$ be the time that vehicle $i$ enters the merging zone $z, z=1,2$, and $t_i^{f}$ be the time that vehicle $i$ exits the last merging zone along its route. 
Each vehicle is modeled as a double integrator,
\begin{equation}%
\begin{split}
\dot{p}_{i} =v_{i}(t), \quad \dot{v}_{i} =u_{i}(t),
\label{eq:model2}
\end{split}
\end{equation}
where $p_{i}(t)\in\mathcal{P}_{i}$, $v_{i}(t)\in\mathcal{V}_{i}$, and $u_{i}(t)\in\mathcal{U}_{i}$ denote the position, speed and
acceleration of each vehicle $i$ in the control zone of the corridor. The sets $\mathcal{P}_{i}$, $\mathcal{V}_{i}$, and $\mathcal{U}_{i}$, $i\in\mathcal{N}(t),$ are complete and totally bounded subsets of $\mathbb{R}$.
Let $x_{i}(t)=\left[p_{i}(t) ~ v_{i}(t)\right]  ^{T}$ denote the state of each vehicle $i$, with initial value $x_{i}^{0}=x_i(t_i^0)=\left[p_{i}^{0} ~ v_{i}^{0}\right]  ^{T}$, where $p_{i}^{0}= p_{i}(t_{i}^{0})=0$ at the entry of the corridor, taking values in $\mathcal{X}_{i}%
=\mathcal{P}_{i}\times\mathcal{V}_{i}$. The state space $\mathcal{X}_{i}$ for each vehicle $i$ is closed with respect to the induced topology on $\mathcal{P}_{i}\times \mathcal{V}_{i}$ and thus, it is compact.

We need to ensure that for any initial state $(t_i^0, x_i^0)$ and every admissible control $u_i(t)$, the system \eqref{eq:model2} has a unique solution $x_i(t)$ on some interval $[t_i^0, t_i^f]$. The following observations from \eqref{eq:model2} satisfy some regularity conditions required both on the model and admissible controls $u_i(t)$ to guarantee local existence and uniqueness of solutions for \eqref{eq:model2}: a) The model is continuous in $u$ and continuously differentiable in the state $x$, b) The first derivative of the model in $x$ is continuous in $u$, and c) The admissible control $u_i(t)$ is continuous with respect to $t$.

To ensure that the control input and vehicle speed are within a
given admissible range, the following constraints are imposed,
\begin{equation}%
\begin{split}
u_{i,min}   &\leq u_{i}(t)\leq u_{i,max},\quad\text{and}\\
0   \leq v_{min}&\leq v_{i}(t)\leq v_{max},\quad\forall t\in\lbrack t_{i}%
^{0},t_{i}^{f}],
\end{split}
\label{eq:speed_accel constraints}%
\end{equation}
where $u_{i,min}$, $u_{i,max}$ are the minimum and maximum
acceleration for each vehicle $i\in\mathcal{N}(t)$, and $v_{min}$, $v_{max}$ are the minimum and maximum speed limits respectively. 

To characterize the physical location of the CAV $i-1\in \mathcal{Q}^z(t)$ inside the control zone, three subsets $\mathcal{L}_i^{z}(t),\mathcal{O}_i^{z}(t)$ and $\mathcal{C}_i^{z}(t)$ of $\mathcal{Q}^z(t)$ with respect to CAV $i$ are defined as follows.

\textit{Definition 2:}
\label{def: subsets}
1) $\mathcal{L}_i^{z}(t)$ contains all vehicles traveling in the same direction and same lane as vehicle $i$ with a potential of rear-end collision, e.g., $\mathcal{L}_5^1(t)$ contains CAV $\# 4$ (Fig. \ref{fig:1}), 2) $\mathcal{O}_i^{z}(t)$ contains all vehicles that travel in the opposite direction as vehicle $i$, and thus no rear-end or lateral collision is possible, e.g., $\mathcal{O}_1^2(t)$ contains CAV $\# 3$ (Fig. \ref{fig:1}), and 3) $\mathcal{C}_i^{z}(t)$ contains all vehicles from different entry points with the possibility of  lateral collision with vehicle $i$, e.g., $\mathcal{C}_8^2(t)$ contains CAV $\# 1$ (Fig. \ref{fig:1}).

To ensure the absence of rear-end collision of two consecutive vehicles traveling on the same lane, the position of the physically immediately preceding CAV  $k\in \mathcal{L}_i^{z}(t)$ should be greater than or equal to the position of the following vehicle plus a predefined safe distance headway $\delta_i(t)$, which is a function of speed $v_i(t)$.
Thus we impose the rear-end safety constraint, 
\begin{equation}
\begin{split}
p_{k}(t)-p_{i}(t) \ge \delta_i(t),~ \forall t\in [t_i^0, t_i^f].
\label{eq:rearend_constraint}
\end{split}
\end{equation}

For each CAV $i\in \mathcal{Q}^z(t)$, the lateral collision is possible within the set $\Gamma_i$,
\begin{equation}
\Gamma_i \overset{\Delta}{=} \{t\,\,|\,t\in  [t_i^{m_z}, t_i^{f}]\}.
\end{equation}
Lateral collision between any two CAVs $i,j\in \mathcal{Q}^z(t)$ can be avoided if the following constraint hold
\begin{equation}\label{eq:lateral_constraint}
\Gamma_i \cap \Gamma_j=\varnothing, ~  \forall t\in [t_i^{m_z}, t_i^f],\,\, i,j\in \mathcal{Q}^z(t).
\end{equation} 


\subsection{Upper Level Vehicle Coordination Problem}
The upper level vehicle coordination problem provides the time $t_i^{m_z^*}$ that each CAV $i\in \mathcal{Q}^z(t)$ enters the merging zone $z, z=1,2$. For $i=1$, due to the absence of any prior CAV inside the control zone, the safety constraints \eqref{eq:rearend_constraint} and \eqref{eq:lateral_constraint} are  not active. This leads to the trivial solution $v_1^*(t)=v_i^0, \quad \forall t\in [t_i^0, t_i^{m_z}]$ and $t_i^{m_z^*} =  t_i^0+\frac{L_z}{v_i(t_i^0)}$. For the rest of the CAVs $i\in \mathcal{Q}^z(t)$, we seek to maximize the traffic throughput by minimizing the inter-vehicle gaps according the following optimization scheme,

\begin{equation}
\begin{split}
\label{eq:throughput_max}
\underset{t(2:M^z(t))}{\text{min}} \sum_{z=1}^{2}\sum_{i=2}^{M^z(t)}(t_i^{m_z} -t_{i-1}^{m_z})=\\
\underset{t(M^z(t))}{\text{min}}\sum_{z=1}^{2}(t_{M^z(t)}^{m_z}-t_1^{m_z}) \\ 
\text{subject to : } \, \eqref{eq:MS}, \eqref{eq:speed_accel constraints}, \eqref{eq:rearend_constraint}, \eqref{eq:lateral_constraint}. 
\end{split}
\end{equation} 

The solution of the optimal control problem \eqref{eq:throughput_max} recursively yields a feasible merging time $t_i^{m_z^*}$ for each vehicle $i$ to cross the merging zone $z$ satisfying condition \eqref{eq:MS}

\textit{Theorem 1:}
If the state and control constraints in \eqref{eq:speed_accel constraints} are inactive, the solution $\mathcal{T}^*=\{t_2^{m_z^*},...,t_{M^z(t)}^{m_z^*} \}$ of \eqref{eq:throughput_max} can be obtained through the following recursive structure over $i=2,...,M^z(t)$ for each $z, z=1,2$,

\begin{equation}\label{eq:veh-coordination}
t_i^{m_z^*} =\left\{
\begin{array}
[c]{ll}%
t_{k}^{m_z^*}+\frac{\delta_i(t)}{v_{k}(t_{k}^{m_z^*})},\text{if} \,\,i-1\in \mathcal{L}_i^{z}(t) ,\\
\text{max}\bigg\{t_{i-1}^{m_z^*},t_{k}^{m_z^*} +\frac{\delta_i(t)}{v_{k}(t_{k}^{m_z^*})}\bigg\},\\ \quad\quad\quad \text{if} \,\, i-1\in \mathcal{O}_i^{z}(t)\\
\text{max}\bigg\{t_{i-1}^{m_z^*} + \frac{S_z}{v_{i-1}(t_{i-1}^{m_z^*})},\\ \quad \quad t_{k}^{m_z^*} +\frac{\delta_i(t)}{v_{k}(t_{k}^{m_z^*})}\bigg\},\text{if} \,\,i-1\in \mathcal{C}_i^{z}(t)\\
\end{array}
\right.
\end{equation}
Due to space limitations, the proof is not presented here. However, the structure and the steps of the proof are similar to Theorem 1 in \cite{Malikopoulos2017}.
If condition \eqref{eq:MS_check} is violated, two special cases arise. Based on the following proposition, CAV $i$ can either follow $i-1$ or reach the merging zone before $i-1$.
 
\textit{Proposition 1:} If there exists a CAV $j\in \mathcal{C}_i^{z}(t):\left|t_i^{m_z}-t_j^{m_ z^*}\right| < \rho_i$ or $k\in \mathcal{L}_i^{z}(t)$, the order of the CAVs in $\mathcal{Q}^{z}(t)$ is conserved and $t_{i}^{m_z^*}$ is calculated by \eqref{eq:veh-coordination}. If there is no CAV $k\in \mathcal{L}_i^{z}(t)$ and $\left|t_i^{m_z}-t_j^{m_ z^*}\right| \ge \rho_i, \forall j \in \mathcal{C}_i^{z}(t)$, then the order of the CAVs in $\mathcal{Q}_i^{z}(t)$ is updated and $t_i^{m_z^*} = t_i^{0}+\frac{L_z}{v_i(t_i^0)} $ Here, $\rho_i$ is a predefined safe time headway.

\begin{proof}
Part 1: If there exists a CAV $ k\in \mathcal{L}_i^{z}(t)$, or $j\in \mathcal{C}_i^{z}(t):\left|t_i^{m_z}-t_j^{m_ z^*}\right| < \rho_i$, CAV $i$ cannot have a collision free trajectory. This implies $i$ cannot travel with its initially calculated merging time $t_i^{m_z}$ and has to merge after CAV $i-1$. To satisfy \eqref{eq:MS}, the merging sequence is conserved, and $t_i^{m_z^*}$ is calculated by \eqref{eq:veh-coordination} which minimizes \eqref{eq:throughput_max}. 

Part 2: With the absence of $k\in \mathcal{L}_i^{z}(t)$, if  $\left|t_i^{m_z}-t_j^{m_z^*}\right| \ge \rho_i, \forall j \in \mathcal{C}_i^{z}(t)$, vehicle $i$ can have a collision-free trajectory and maintain its initial velocity such that $t_i^{m_z^*}=t_i^{m_z}=t_i^{0}+\frac{L_z}{v_i(t_i^0)}$. In this case, the merging sequence is updated so that \eqref{eq:MS} is not violated.
\end{proof}

\subsection{Lower Level Energy Minimization Problem}
For each vehicle $i\in\mathcal{N}(t)$, we formulate the decentralized optimial control problem that minimizes the cost function $J_{i}(u(t))$ in $[t_i^0, t_{i}^{m_z}]$,
\begin{gather}\label{eq:decentral_general}
\min_{u_i\in U_i}J_{i}(u(t))=  \int_{t^0_i}^{t_{i}^{m_z}} C_i(u_i(t))~dt,\\ 
\text{subject to}%
:\eqref{eq:model2},\eqref{eq:speed_accel constraints},\text{
}p_{i}(t_{i}^{0})=0\text{, }p_{i}(t_{i}^{m_z})=L_z,\nonumber\\
\text{and given }t_{i}^{0}\text{, }v_{i}^{0}\text{, }t_{i}^{m_z}.\nonumber
\end{gather}
Here, $C_i(u_i(t))$ is monotonically increasing function of the control input $u_i(t)$ and can be viewed as a measure of energy. When $C_i(u_i(t))$ is considered as the $L^2$-norm of the control input, i.e. $C_i(u_i(t))=\frac{1}{2}u_i^2(t)$, the transient engine operation can be minimized, which eventually represents the minimization of fuel consumption \cite{Malikopoulos2010a}.  Note that, the safety constraints are not included in \eqref{eq:decentral_general}. The lateral collision constraint \eqref{eq:lateral_constraint} is implicitly included by solving the upper-level vehicle coordination problem. The rear-end collision constraint \eqref{eq:rearend_constraint} can be avoided under proper initial conditions $[t_i^0,v_i^0(t)]$ as described in \cite{Malikopoulos2017}.

In what follows, we provide the closed-form solution of the optimal control problem formulated in \eqref{eq:decentral_general} for each vehicle $i\in\mathcal{N}(t)$.

\section{Analytical Solution of the Lower Level Energy Minimization Problem} \label{sec:3} 
The solution of the constrained problem has been addressed in \cite{Malikopoulos2017}, and it requires the constrained and unconstrained arcs of the state and control input to be pieced together to satisfy the Euler-Lagrange equations and necessary condition of optimality. Due to the page limitations, we provide the general formulation and include only the solution of the unconstrained case here. From \eqref{eq:decentral_general}, the state equations \eqref{eq:model2} and the constraints \eqref{eq:speed_accel constraints}, for
each vehicle $i\in\mathcal{N}(t)$ the Hamiltonian function with the state and control adjoined is
\begin{gather}
H_{i}\big(t, x(t),u(t)\big) = \frac{1}{2} u^{2}_{i} + \lambda^{p}_{i} \cdot
v_{i} + \lambda^{v}_{i} \cdot u_{i} %
\nonumber\\
+ \mu^{a}_{i} \cdot(u_{i} - u_{max})+ \mu^{b}_{i} \cdot(u_{min} - u_{i}) + \mu^{c}_{i} \cdot(v_{i} - v_{max}) \nonumber\\
+ \mu^{d}_{i} \cdot(v_{min} - v_{i}), \quad \forall i \in\mathcal{N}(t),\label{eq:16b}%
\end{gather}
where $\lambda^{p}_{i}$ and $\lambda^{v}_{i}$ are the co-state components, and
$\mu_i^{a},\mu_i^{b},\mu_i^{c}$ and $\mu_i^{d}$ are the Lagrange multipliers.

%

\subsection{Analytical Solution without Interior Constraints}
If the inequality state and control constraints \eqref{eq:speed_accel constraints} are not active, we have $\mu^{a}%
_{i} = \mu^{b}_{i}= \mu^{c}_{i}=\mu^{d}_{i}=0$. Applying the necessary condition, the optimal control can be given 
\begin{equation}
u_i(t) + \lambda^{v}_{i}= 0, \quad i \in\mathcal{N}(t). \label{eq:17}
\end{equation}
From Euler-Lagrange equations, we have $\lambda^{p}_{i}(t) = a_{i}$, and $\lambda^{v}_{i}(t) = -\big(a_{i}\cdot t + b_{i}\big)$. 
The coefficients $a_{i}$ and $b_{i}$ are constants of integration corresponding to each vehicle $i$. From \eqref{eq:17} the optimal control input (acceleration/deceleration) as a function of time, and the corresponding state trajectories are given by

\begin{gather}
u^{*}_{i}(t) = a_{i} \cdot t + b_{i}, ~ \forall t \ge t^{0}_{i}. \label{eq:20}\\
v^{*}_{i}(t) = \frac{1}{2} a_{i} \cdot t^2 + b_{i} \cdot t + c_{i}, ~ \forall t \ge t^{0}_{i}\label{eq:21}\\
p^{*}_{i}(t) = \frac{1}{6}  a_{i} \cdot t^3 +\frac{1}{2} b_{i} \cdot t^2 + c_{i}\cdot t + d_{i}, ~ \forall t \ge t^{0}_{i}, \label{eq:22}%
\end{gather}
where $c_{i}$ and $d_{i}$ are constants of integration. The constants of integration $a_i$, $b_{i}$, $c_{i}$, and $d_{i}$ can be computed once at time $t^{0}_{i}$ using the initial and final conditions, and the values of the one of terminal transversality condition, i.e., $\lambda^{v}_{i}(t_i^{m_z})=0$.

\subsection{Analytical Solution with Interior Constraints}
In the case that the path of vehicle $i$ consists of more than one merging zone, for example, the eastbound CAV $i$ enters from the left and travels through merging zone \#1 and \#2 in Fig. \ref{fig:1} between the time $t_i^0$ that the vehicle enters the control zone and the time $t_i^f$ that the vehicle exits the merging zone \#2, vehicle $i$ has to travel across the intermediate merging zone \#1 at the designated time $t_i^{m_1}$. Therefore, we need to impose an additional interior boundary condition \cite{Malikopoulos2018d}
\begin{gather}
p_i(t_i^{m_1}) = L_1.
\label{eq:interior1}%
\end{gather} 
If a speed constraint $v_1$ is imposed as an interior boundary condition, then
\begin{gather}
v_i(t_i^{m_1}) = v_1.
\label{eq:interior1_v}%
\end{gather} 
Let $t_i^{m_1-}$ and $t_i^{m_1+}$ represents the time just before and after the jump conditions. Then
\begin{gather}
\lambda_i^p(t_i^{m_1-}) = \lambda_i^p(t_i^{m_1+}) + \pi_0,
\label{eq:interior2}\\
\lambda_i^v(t_i^{m_1-}) = \lambda_i^v(t_i^{m_1+}) + \pi_1,
\label{eq:interior3}\\
H^- = H^+ - \pi_0 \cdot v_i(t_i^{m_1}) - \pi_1 \cdot u_i(t_i^{m_1}).
\label{eq:interior4}%
\end{gather}
where $\pi_0$ and $\pi_1$ are  constant Lagrange multipliers, determined so that \eqref{eq:interior2} and \eqref{eq:interior3} are satisfied. Equations \eqref{eq:interior2}- \eqref{eq:interior4} imply discontinuities in the position and speed co-states and the Hamiltonian at  $t_i^{m_1}$. The two arcs, i.e., equations before and after $t_i^{m_1}$, are pieced together to solve the problem with 9 or 10 unknowns [if  \eqref{eq:interior1_v} is also imposed]  including the constants of integration, $\pi_0$ and/or $\pi_1$, and the corresponding equations: the initial conditions, i.e., $v_i(t_i^0)$ and $p_i(t_i^0)$, the interior conditions as defined in (\ref{eq:interior1}), [and/or (\ref{eq:interior1_v})] 
the final conditions, i.e., $\lambda_i(t_i^{m_z}), p_i(t_i^{m_z})$, and the junction point defined in (\ref{eq:interior2}) [and/or (\ref{eq:interior3})].

To derive online the optimal control for each vehicle $i$, we need to calculate the constants of integration at time $t_i^0$, so that the controller yields the optimal control online for each vehicle $i$. We form the following system of nine equations, namely
\begin{gather} \small
\label{eq:matrix}%
\noindent
\renewcommand\arraystretch{1.6} 
\resizebox{\linewidth}{!}{%
   $ \begin{pmatrix}
    \frac{(t_i^0)^2}{2} & (t_i^0) & 1 & 0 & 0 & 0 & 0 & 0 & 0\\
	\frac{(t_i^0)^3}{6} & \frac{(t_i^0)^2}{2} & t_i^0 & 1 & 0 & 0 & 0 & 0 & 0\\
	\frac{(t_i^{m_1})^3}{6} & \frac{(t_i^{m_1})^2}{2} & t_i^{m_1} & 1 & 0 & 0 & 0 & 0 & 0\\
	\frac{(t_i^{m_1})^2}{2} & t_i^{m_1} & 1 & 0 & -\frac{(t_i^{m_1})^2}{2} & -t_i^{m_1} & -1 & 0 & 0\\
	0 & 0 & 0 & 0 & \frac{(t_i^{m_2})^3}{6} & \frac{(t_i^{m_2})^2}{2} & (t_i^{m_2}) & 1 & 0\\
	0 & 0 & 0 & 0 & -t_i^{m_2} & -1 & 0 & 0 & 0\\
	0 & 0 & 0 & 0 & \frac{(t_i^{m_1})^3}{6} & \frac{(t_i^{m_1})^2}{2} & t_i^{m_1} & 1 & 0\\
	t_i^{m_1} & 1 & 0 & 0 & -t_i^{m_1} & -1 & 0 & 0 & 0\\
	1 & 0 & 0 & 0 & -1 & 0 & 0 & 0 & -1
    \end{pmatrix}  $
    } 
    \nonumber\\
	\cdot   
    \begin{pmatrix}
    a_{i}\\
	b_{i}\\
	c_{i}\\
	d_{i}\\
	g_{i}\\
	h_{i}\\
	q_{i}\\
	w_{i}\\
	\pi_0
    \end{pmatrix}
    =
    \begin{pmatrix}
    v_{i}(t_i^0)\\
	p_{i}(t_i^0)\\
	p_{i}(t_i^{m_1})\\
	0\\
	p_{i}(t_i^{m_2})\\
	\lambda^{v}_{i}(t_i^{m_2})\\
	p_{i}(t_i^{m_1})\\
	0\\
	0
    \end{pmatrix}, \forall t \ge t^{0}_{i}.
\end{gather}
where $a_i, b_i, c_i, d_i$ are the constants of integration for the first arc, and $g_i,h_i,q_i,w_i$ are the constants of integration for the second arc.

\section{Simulation Results}

\begin{figure} [ht]
	\centering
	\includegraphics[scale=0.4]{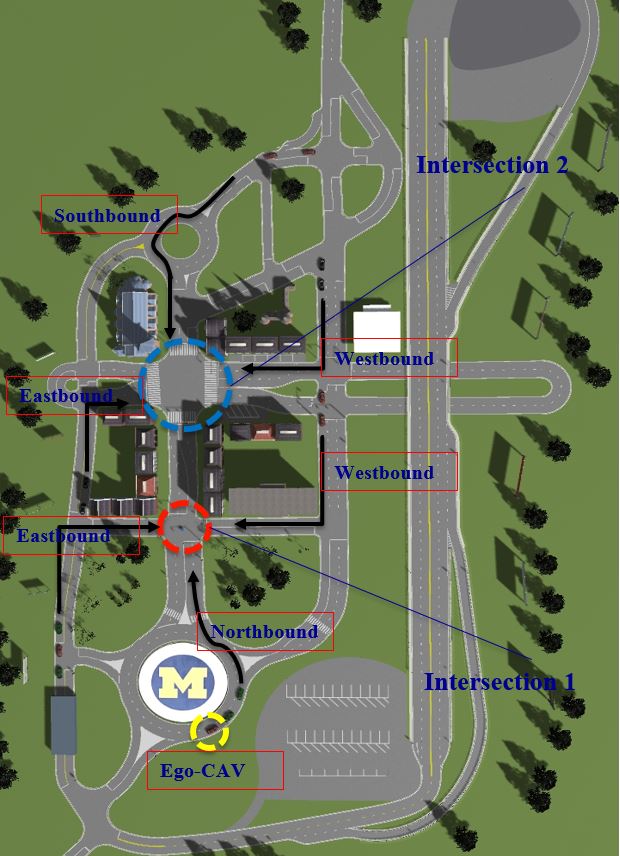}   
	\caption{Illustration of the corridor in Mcity.} 
	\label{fig:corridor}
\end{figure}

To evaluate and validate the effectiveness of the proposed approach, we conducted  computational studies using the commercial software platforms of TASS International PreScan in conjuction with Mathworks MATLAB and Mathworks Simulink. We considered a corridor with two adjacent intersections (intersection-1 and intersection-2) in Mcity (Fig. \ref{fig:corridor}), a 32 acre vehicle testing facility. The dimensions of the conflict zones are 18 $m$ $\times$ 12 $m$ for intersection-1 and 34 $m$ $\times$ 28 $m$ for intersection-2. We select the length of the control zone to be $100$ $m$ measured from the entry of each intersection. Six different routes have been designed for the scenario in Fig. \ref{fig:corridor}  with 14 CAVs: 1) two eastbound routes with 5 CAVs, 2) two westbound routes with 4 CAVs, 3) one southbound route with 2 CAVs, and finally 4) one northbound route with 3 CAVs. Note that, east and westbound vehicles travel through only one intersection in their path. The routes and CAV positions were designed in such a way that the trajectories of the CAVs present a worst-case collision scenario. To analyze the individual performance of the simulated CAVs, we considered a northbound vehicle (ego-CAV) as the test vehicle (see Fig. \ref{fig:corridor}).

\begin{figure} [ht]
	\centering
	\includegraphics[width=0.5\textwidth]{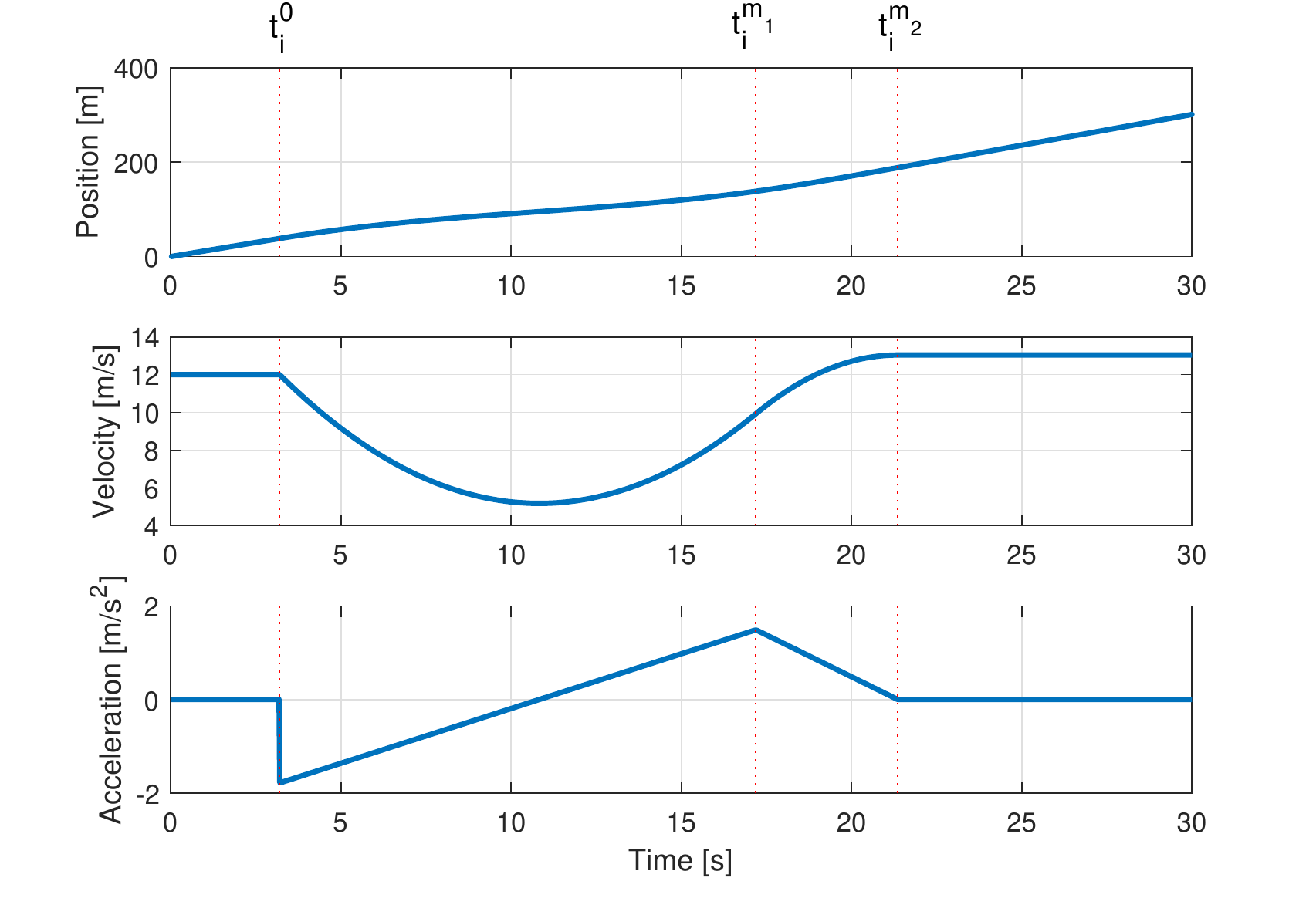}  
	\caption{Optimal vehicle trajectory traveling through two intersections.} 
	\label{fig:prof_MA1}
\end{figure}

\begin{figure} [ht]
	\centering
	\includegraphics[width=0.53\textwidth]{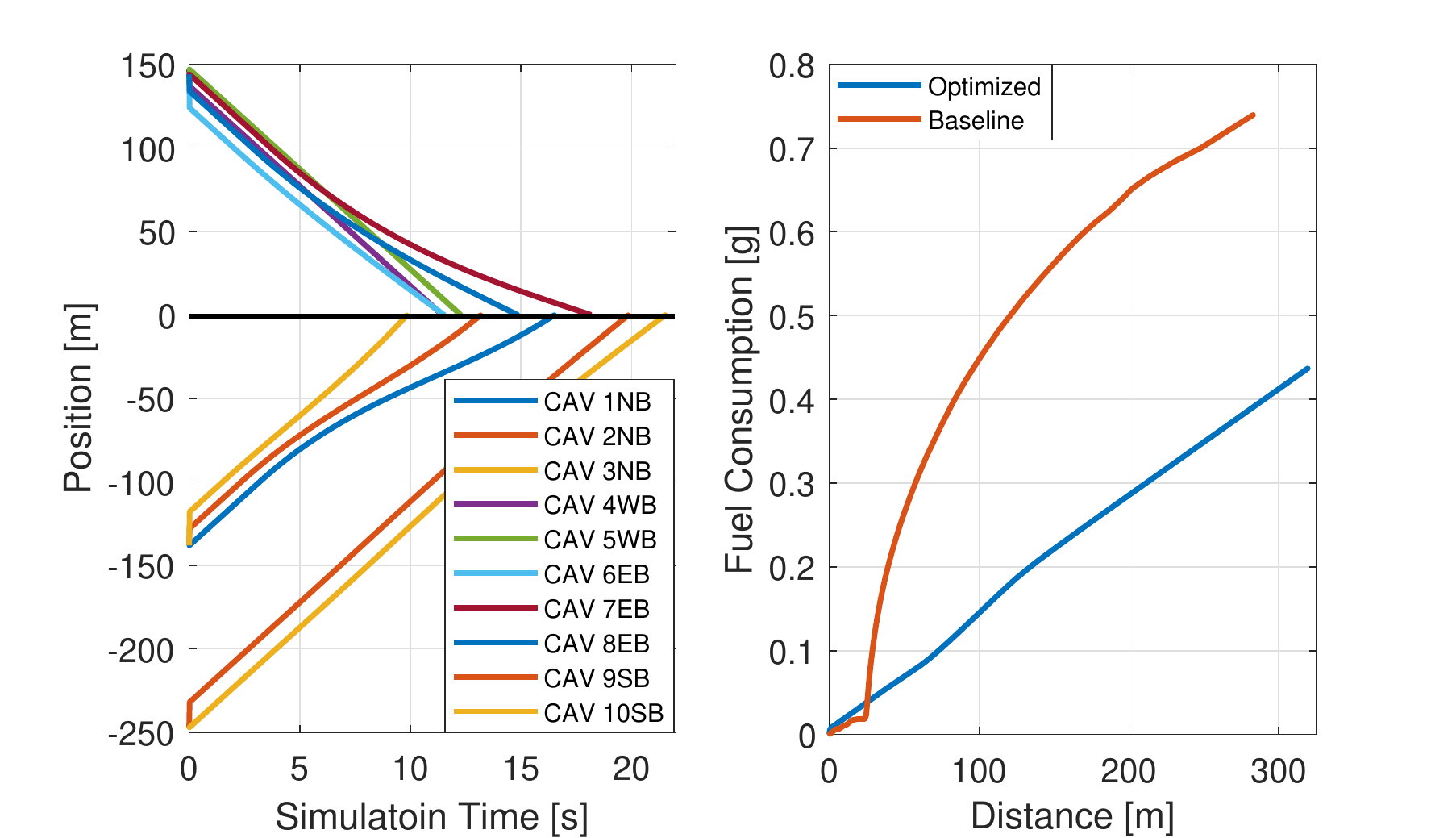}   
	\caption{Collision-free optimal trajectory of 10 CAVs approaching towards intersection-1 (left) and cumulative fuel consumption of the optimized and baseline scenarios (right).} 
	\label{fig:int-FC}
\end{figure}

The merging scenario at intersection-1 with 10 incoming CAVs is depicted in Fig. \ref{fig:int-FC} (left). Note that, CAV$\#4$ heading westbound and CAV$\#7$ heading eastbound are allowed to enter the conflict zone at the same time since their routes have opposite direction and thus are non-conflicting. The simulation results for the decentralized optimization problem for the ego-CAV are depicted in Fig. \ref{fig:prof_MA1}. We observe that, the optimal control takes hold of the ego-CAV at the entry of the control zone denoted at time $t_i^0$ and leads it optimally through intermediate collision points of intersection-1 at time $t_i^{m_1}$ and $t_i^{m_2}$.

To compare the performance of the proposed optimal solution, we construct a baseline scenario with fixed time signalized intersections with switching time of 10 seconds. The vehicles were governed by the $Gipps$ car-following model \cite{Gipps1981}. To quantify the effect of optimal vehicle coordination on fuel consumption, a polynomial meta-model proposed in \cite{Kamal2013a} was used. A comparison of fuel consumption for the ego-CAV between the baseline and optimized scenarios is shown in Fig. \ref{fig:int-FC} (right). We observe $40.9\%$ improvement in fuel efficiency for the  ego-CAV under the constructed baseline scenario.

\begin{figure} [ht]
	\centering
	\includegraphics[width=0.5\textwidth]{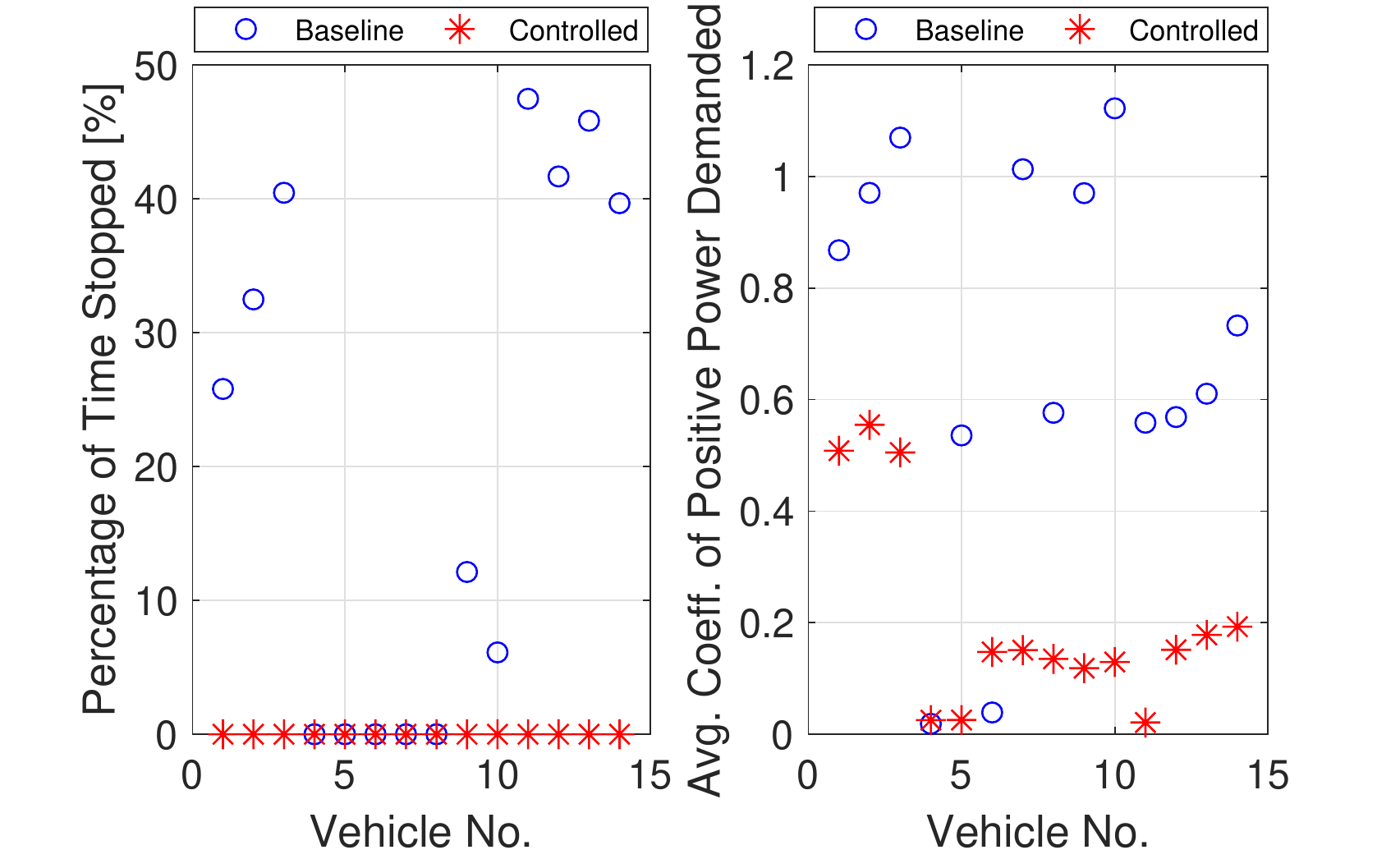}   
	\caption{Stoppage time (left) and average coefficient of power demanded (right) of the optimized and baseline scenarios.} 
	\label{doubleInter_coeff_power_demand}
\end{figure}

\begin{figure} [ht]
	\centering
	\includegraphics[width=0.4\textwidth]{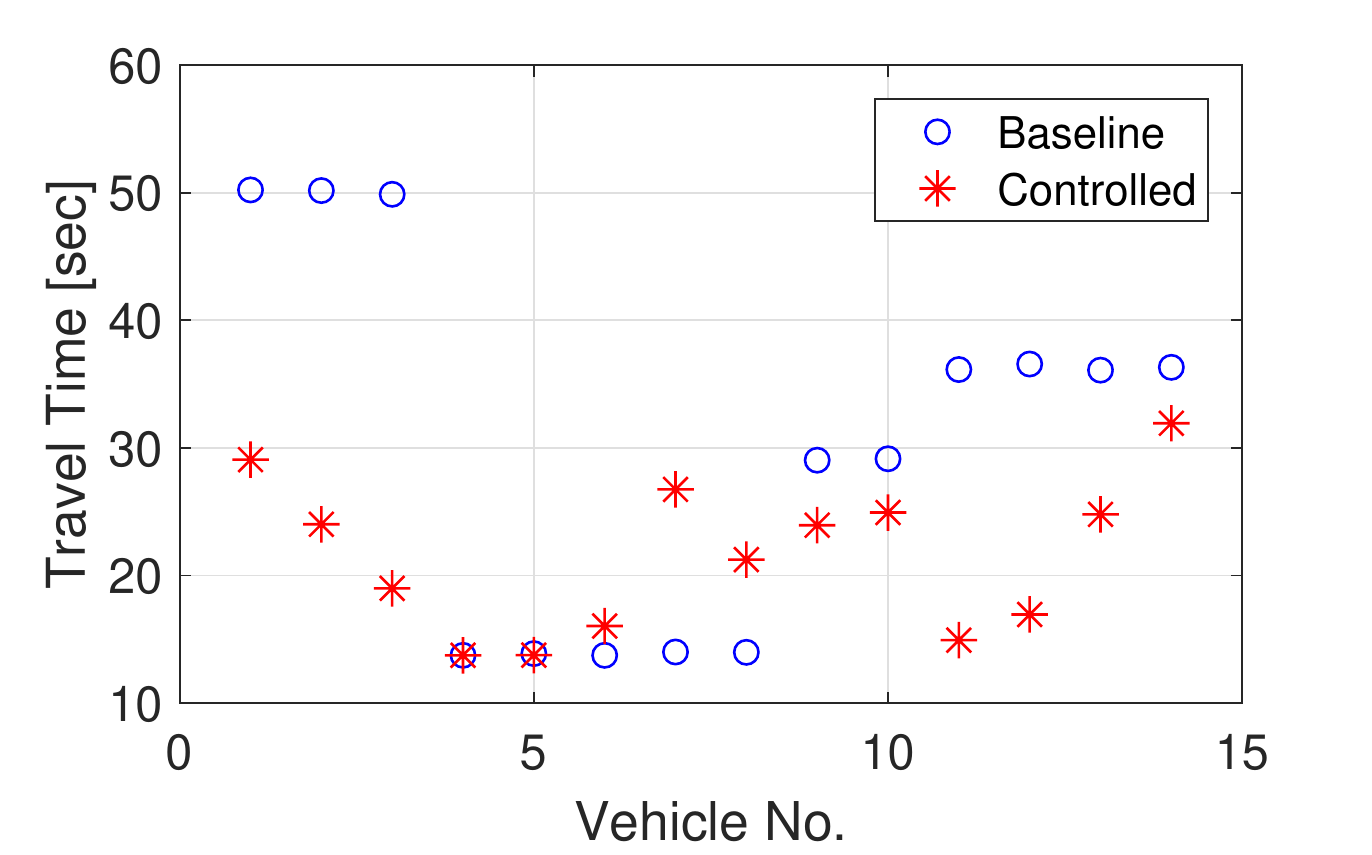}   
	\caption{Fleet travel time of the optimized and baseline scenarios.} 
	\label{doubleInter_travel_time}
\end{figure}

The proposed framework alleviates stop-and-go driving, and thus, minimizes associated transient engine operation in this corridor, yielding improvements in fuel consumption. To quantitatively investigate this observation, the individual drive cycles of the vehicles in the baseline and optimized scenarios are analyzed using three metrics \cite{Malikopoulos2018e,Malikopoulos2013}: (1) total travel time, (2) \textit{stop factor}, and (3) \textit{average coefficient of power demanded}. The stop factor provides a convenient indication of idle engine operation over a driving cycle. The coefficient of power demanded provides an indication of the transient engine operation since it is proportional to power demanded by the driver. Total stoppage time in the drive cycle, shown in Fig. \ref{doubleInter_coeff_power_demand} (left), was eliminated for every vehicle in the fleet. The coefficient of power demand, shown in Fig. \ref{doubleInter_coeff_power_demand} (right), only considers vehicle power demanded under both positive acceleration and velocity events and was able to be reduced by 40.8\% across the fleet of 14 vehicles. The total travel time (Fig. \ref{doubleInter_travel_time}) for all 14 vehicles was improved by 13.2\% with the proposed framework compared to the baseline scenario.

\section{Concluding Remarks}
In this paper, we addressed the problem of coordinating CAVs at two signal-free adjacent intersections by formulating a decentralized optimal control problem. We presented a closed-form analytical solution that considers interior boundary conditions and provides optimal fuel-efficient and collision-free trajectories to the CAVs for their predetermined routes. The proposed decentralized framework exhibits significant improvement in terms of fuel efficiency, average power demand and average travel time when compared to the baseline scenario. Ongoing  efforts consider the complete solution that includes state and control constraints with left/right turns and lane changes. Future research should focus on vehicle coordination under mixed traffic environment where the interaction between human-driving vehicles and CAVs is taken into consideration.

\section{Acknowledgments}
The authors would like to thank Ioannis V. Chremos for the discussions on the Theorem and Proposition.
\bibliographystyle{IEEEtran}
\bibliography{itsc_corridor_ref}

%

%

\end{document}